\documentclass[11pt,leqno]{article}
\usepackage{amsmath}
\usepackage{amssymb}
\usepackage{dsfont}
\usepackage{multicol}
\usepackage{latexsym}
\usepackage[usenames]{color}
\usepackage[colorlinks=true,
linkcolor=webgreen, filecolor=webbrown,
citecolor=webgreen]{hyperref}
\definecolor{webgreen}{rgb}{0,.5,0}
\definecolor{webbrown}{rgb}{.6,0,0}

\usepackage{graphics,amsmath,amssymb}
\usepackage{amsthm}
\usepackage{latexsym}

\def\N{{\mathds{N}}}
\def\Z{{\mathds{Z}}}

\def\d{{\rm{d}}}

\newtheorem{thm}{Theorem}

\begin{document}

\title{\bf A modified M\"obius $\mu$-function}
\author{{\sc Rasa Steuding}, {\sc J\"orn Steuding}, {\sc L\'aszl\'o T\'oth}}
\date{}
\maketitle

\centerline{Rend. Circ. Mat. Palermo {\bf 60} (2011), 13--21}
\vskip4mm

\centerline{Rasa \& J\"orn Steuding} \centerline{Department of
Mathematics, W\"urzburg University} \centerline{Am Hubland, 97\,218
W\"urzburg, Germany}
\centerline{\texttt{steuding@mathematik.uni-wuerzburg.de}}

\vskip1mm \centerline{L\'aszl\'o T\'oth} \centerline{Department of
Mathematics, University of P\'ecs} \centerline{7624 P\'ecs,
Ifj\'us\'ag u. 6, Hungary}
\centerline{\texttt{ltoth@gamma.ttk.pte.hu}}

\begin{abstract} We investigate a modified M\"obius $\mu$-function which is related to
an infinite product of shifted Riemann zeta-functions. We prove
conditional and unconditional upper and lower bounds for its
summatory function, and, finally, we discuss relations with
Riemann's hypothesis.
\end{abstract}



\section{Introduction and prehistory}

The classical M\"obius $\mu$-function is defined by $\mu(1)=1,\
\mu(n)=0$ if $n$ has a quadratic divisor $\neq 1$, and
$\mu(n)=(-1)^r$ if $n$ is the product of $r$ distinct primes. It is
easily seen that $\mu(n)$ is multiplicative and appears as
coefficients of the Dirichlet series representation of the
reciprocal of the Riemann zeta-function:
$$
\zeta(s)^{-1}=\prod_{p\,{\tiny\rm{prime}}}\left(1-{1\over p^s}\right)=\sum_{n=1}^\infty{\mu(n)\over n^s},
$$
both representations being valid for $\sigma>1$, where $s=\sigma+it$
with $i:=\sqrt{-1}$ is a complex variable. Riemann's famous open
hypothesis on the non-vanishing of $\zeta(s)$ in the half-plane
$\sigma>{1\over 2}$ is known to be equivalent to the estimate
$$
{\sf M}(x):=\sum_{n\leq x}\mu(n)\ll x^{1/2+\epsilon}
$$
for any positive $\epsilon$. Odlyzko \& te Riele \cite{or} disproved the original
Mertens hypothesis \cite{merte}, that is $\vert {\sf M}(x)\vert< x^{1/2}$, by showing
$$
\liminf_{x\to\infty}{{\sf M}(x)\over x^{1/2}}<-1.009\qquad\mbox{and}\qquad
\limsup_{x\to\infty}{{\sf M}(x)\over x^{1/2}}>1.06\, ;
$$
for more details see Titchmarsh \cite{titch} (incl. the notes to \S 14).

In this note we are concerned with asymptotic properties of a
modified M\"obius function which is defined as the multiplicative
arithmetical function $\mu_{\infty}$  given by
\begin{equation}\label{uno}
\mu_{\infty}(p^\nu)= (-1)^{|B(\nu)|},
\end{equation}
for any prime power $p^{\nu}$, $\nu \in\N$, where $|B(\nu)|$ is the
number of nonzero terms in the binary representation of the integer
$\nu$, i.e., $\nu=\sum_{j\in B(\nu)} 2^j$. Here
$\mu_{\infty}(p)=\mu_{\infty}(p^2)=-1$, $\mu_{\infty}(p^3)=1$,
$\mu_{\infty}(p^4)=-1$, $\mu_{\infty}(p^5)=\mu_{\infty}(p^6)=1$,
$\mu_{\infty}(p^7)=-1$, etc. This arithmetical function was
introduced by Cohen and Hagis \cite{Coh_Hag1993} and it is an
interesting function for several reasons.

First of all, $\mu_{\infty}$ is the inverse of the function constant
$1$ under the infinitary convolution given by
\begin{equation}\label{infty}
(f\times_{\infty} g)(n)= \sum_{d\mid_{\infty} n} f(d)g(n/d),
\end{equation}
where the sum is over the so called infinitary divisors of $n$,
which are defined in the following way. The infinitary divisors of
the integer $n=\prod p^\nu>1$ are number $1$ and the products of
prime power divisors of $n$ of the form $p^{2^j}$, where $j\in
B(\nu)$ with the notation of above. By convention, $1\mid_{\infty}
1$. The term ``infinitary'' is justified by an equivalent definition
given by Cohen \cite{Coh1990}. By a curious property,
$p^b\mid_{\infty} p^a$ holds if and only if the binomial coefficient
$\binom{a}{b}$ is odd.

On the other hand, the function $\mu_{\infty}$ is identical with the function
denoted by $\mu^{**}$, which is the inverse of the function constant
$1$ under the bi-unitary convolution defined by
\begin{equation}\label{zwo}
(f \times_{**} g)(n)= \sum_{\substack{d\mid n\\ (d,n/d)_{**}=1}}f(d)g(n/d),
\end{equation}
where $(a,b)_{**}$ denotes the greatest common unitary divisor of
$a$ and $b$. Recall that $e$ is said to be a unitary divisor of $k$
if $e$ divides $k$ with greatest common divisor $(e,k/e)=1$.  The
sum in (\ref{zwo}) is over the so called bi-unitary divisors of $n$.
The bi-unitary divisors of a prime power $p^a$ ($a\ge 1$) are all
divisors $p^b$ with $b=0,1,2,\ldots,a$, except $p^{a/2}$ for $a$
even. The concept of bi-unitary divisor is due to Suryanarayana
\cite{Sur1972}, while properties of the bi-unitary convolution are
given by Haukkanen \cite{Hau1998}.

Note that both the infinitary and bi-unitary convolutions are
commutative. The infinitary convolution is associative, however the
bi-unitary convolution is not associative. The identity with respect
to both convolutions is the function $\delta$ given by $\delta(1)=1$
and $\delta(n)=0$ for $n>1$. Furthermore, $f$ has an inverse under
each convolutions if and only if $f(1)\ne 0$. If $f$ and $g$ are
multiplicative, then $f\times_{\infty} g$ and $f \times_{**} g$ are
also multiplicative. Moreover, if $f$ is a non-zero multiplicative
function, then their inverses under each convolutions are also
multiplicative, cf. \cite{Coh_Hag1993,Hau1998,Hau2000}.

Besides $\mu_{\infty}$ Cohen and Hagis \cite{Coh_Hag1993} also investigated
the functions $\tau_{\infty}(n)$ and $\sigma_{\infty}(n)$, denoting the
number and the sum of the infinitary divisors of $n$, proving asymptotic
formulae for the summatory functions of $\tau_{\infty}(n)$ and
$\sigma_{\infty}(n)$. Asymptotic formulae for
the corresponding bi-unitary functions $\tau^{**}(n)$
and $\sigma^{**}(n)$ were established in papers
\cite{Sur1972,SurSit1975,SurSub1980}. All of these functions are multiplicative.

These asymptotic formulae may be compared with those involving the classical divisor
function $\tau(n)$ and the sum--of--divisors function $\sigma(n)$.
In this note we shall prove several results concerning the summatory
function of $\mu_{\infty}$ similar to those for the classical M\"obius function.

\section{Main results}

The analytic method has proved to be a rather powerful
approach to study the classical M\"obius $\mu$-function.
We shall mimic this approach and prove first a representation for
the generating Dirichlet series:

\begin{thm}
For $\sigma>1$,
$$
\mathfrak{m}(s):=\sum_{n=1}^\infty {\mu_{\infty}(n)\over n^s}=\prod_{j=0}^\infty \zeta(2^js)^{-1}.
$$
\end{thm}

\noindent {\bf Proof.} Expanding into an Euler product and using (\ref{uno}),
\begin{equation*}
\sum_{n=1}^{\infty} \frac{\mu_{\infty}(n)}{n^s}=  \prod_p \sum_{\nu=0}^{\infty}
\frac{\mu_{\infty}(p^{\nu})}{p^{\nu s}}
= \prod_p \sum_{\nu=0}^{\infty} \frac{(-1)^{|B(\nu)|}}{p^{\nu s}} =
\prod_p \sum_{\nu=0}^{\infty} \prod_{j\in B(\nu)} \frac{-1}{p^{2^j s}}
\end{equation*}
\begin{equation*}
= \prod_p \prod_{j=0}^{\infty} \left(1-\frac1{p^{2^js}} \right)
= \prod_{j=0}^{\infty} \zeta(2^js)^{-1}.
\end{equation*}
The theorem is proved. $\bullet$

We are mainly interested in the asymptotic behaviour of the
summatory function of the modified M\"obius $\mu$-function,
$$
\mathfrak{M}(x):=\sum_{n\leq x}\mu_{\infty}(n),
$$
as $x\to\infty$. Our first theorem is unconditional:

\begin{thm}\label{eins}
There exists a positive constant $c$ such that
\begin{equation}\label{BigO}
\mathfrak{M}(x)=O\left(x\exp\big(-c(\log x)^{3/5}(\log\log x)^{-1/5}\big)\right),
\end{equation}
and, for any positive $\epsilon$,
$$
\mathfrak{M}(x)=\Omega(x^{\beta-\epsilon}),
$$
where $\beta$ is the supremum over all real parts of $\zeta$-zeros
(hence ${1\over 2}\leq \beta\leq 1$).
\end{thm}

Here $f=\Omega(g)$ denotes the negation of $f=o(g)$. The proof of
the first assertion follows along the lines of the proof of the
prime number theorem; the second statement is rather similar to the
so-called Mertens conjecture or how the size of the summatory
function of the M\"obius $\mu$-function is related to the zeros of
the zeta-function.

\vskip1mm \noindent {\bf Proof.} We start with the big-Oh estimate.
By Perron's formula, for $c>1$,
$$
\mathfrak{M}(x)={1\over 2\pi i}\int_{c-iT}^{c+iT}\mathfrak{m}(s){x^s\over s}\d s+{\mathcal E},
$$
where
\begin{equation}\label{e}
{\mathcal E}=O\Big({x^c\over T}\sum_{n=1}^\infty
{\vert\mu_{\infty}(n)\vert\over n^c}+ {x\log x\over T}\Big).
\end{equation}

We shall move the path of integration to the left. Korobov
\cite{koro} and Vinogradov \cite{vino} (independently) proved
$$
\zeta(s)\neq 0\qquad\mbox{in}\quad \sigma\geq 1- C(\log \vert
t\vert+3)^{-2/3}(\log\log (\vert t\vert+3))^{-1/3},
$$
where $C$ is some positive absolute constant; moreover, in the same region the estimate
\begin{equation}\label{schnee}
\zeta(\sigma+it)^{-1}\ll (\log \vert t\vert+3)^{2/3}(\log\log (\vert t\vert+3))^{1/3}
\end{equation}
holds. The first complete proof due to Richert appeared in Walfisz
\cite{walfisch} (see also \cite{ivic}, \S 12). Denote the
rectangular contour with vertices $c\pm iT, 1-\Delta\pm iT$ by
${\mathcal C}$, where $\Delta:={C\over 2}(\log \vert
T\vert+3)^{-2/3}(\log\log (\vert T\vert+3))^{-1/3}$. Then there are
no $\zeta$-zeros on or in the interior of ${\mathcal C}$, and we
deduce from Cauchy's theorem
$$
\mathfrak{M}(x)={1\over 2\pi i}\left\{\int_{c-iT}^{1-\Delta-iT}+
m\int_{1-\Delta-iT}^{1-\Delta+iT}+\int_{1-\Delta+iT}^{c+iT}\right\}\mathfrak{m}(s){x^s\over
s}\d s+{\mathcal E}.
$$

In order to bound the appearing integrals we note for $\sigma>1$
\begin{equation}\label{triv}
{\zeta(2\sigma)\over \zeta(\sigma)}\leq \vert\zeta(s)\vert\leq \zeta(\sigma);
\end{equation}
these inequalities follow factorwise from the Euler product
representation of the zeta-function. For any non-negative integer
$J$ let
\begin{equation}\label{mn}
\mathfrak{m}(s)=N_J(s)\prod_{0\leq j<J}\zeta(2^js)^{-1}\ , \quad
\mbox{where}\quad N_J(s):=\prod_{j\geq J}\zeta(2^js)^{-1}.
\end{equation}

In view of (\ref{triv}), for $\sigma>2^{-J}$,
\begin{equation}\label{tannenbaum}
\vert N_J(s)\vert\leq \prod_{j\geq J}^\infty {\zeta(2^j\sigma)\over \zeta(2^{j+1}\sigma)}=\zeta(2^J\sigma),
\end{equation}
hence the function $N_J(s)$ is bounded for $\sigma >2^{-J}$. Thus we find via (\ref{schnee})
$$
\int_{1-\Delta\pm iT}^{c\pm iT}\mathfrak{m}(s){x^s\over s}\d s=
\int_{1-\Delta\pm iT}^{c\pm iT}\zeta(s)^{-1}N_1(s){x^s\over s}\d
s\ll {x^c\over T\Delta},
$$
and, similarly,
$$
\int_{1-\Delta-iT}^{1-\Delta+iT}\mathfrak{m}(s){x^s\over s}\d s\ll
x^{1-\Delta}\int_0^T\vert\zeta(\sigma+it)\vert^{-1}{\d t\over
1+\vert t\vert}\ll x^{1-\Delta}{\log T\over \Delta}.
$$

Collecting together, we arrive at
\begin{equation}\label{star}
\mathfrak{M}(x)\ll {x^c\over T\Delta}+x^{1-\Delta}{\log T\over \Delta}+{\mathcal E}.
\end{equation}
With $c=1+(\log x)^{-1}$ we have $\sum
\mu_{\infty}(n)n^{-c}\ll\zeta(c)\ll  \log x$ in the estimate for
${\mathcal E}$; choosing $T$ such that $T\log T=x^\Delta$, we obtain
(\ref{BigO}).

Now we prove the big-Omega result. For $\sigma>1$, we find by partial summation
\begin{equation}\label{zuz}
\sum_{n>x}{\mu_{\infty}(n)\over n^s}=-{\mathfrak{M}(x)\over
x^s}+s\int_x^\infty  \mathfrak{M}(u)u^{-s-1}\d u.
\end{equation}

Assuming $\mathfrak{M}(x)=o(x^{\alpha})$ for some positive $\alpha$,
as $x\to\infty$, the right hand-side converges for $\sigma>\alpha$.
Hence, $\mathfrak{m}(s)$ has a convergent Dirichlet series
representation for $\sigma>\alpha$, and thus defines an analytic
function in this half-plane. For $\alpha<\beta$ this contradicts the
poles of $\mathfrak{m}(s)$ at the nontrivial zeros of $\zeta(s)$ on
the critical line (see \cite{titch}, \S 10.2). The theorem is
proved. $\bullet$

An alternative proof is based on the function
$$
F(s):=\sum_{n=1}^\infty {\mu(n)-\mu_{\infty}(n)\over
n^s}=\zeta(s)^{-1} \left(1-\prod_{j=1}^\infty
\zeta(2^js)^{-1}\right);
$$
the analytic behaviour of $F(s)$ implies (via Perron's formula) that
its  summatory function $\sum_{n\leq x}(\mu(n)-\mu_{\infty}(n))={\sf
M}(x)-\mathfrak{M}(x)$ is small, from which one deduces the
estimates of Theorem \ref{eins} by corresponding ones for ${\sf
M}(x)$.

Next, we shall prove an explicit formula for $\mathfrak{M}(x)$
subject to  the truth of the Riemann hypothesis.

\begin{thm}
The Riemann hypothesis is true if and only if $\mathfrak{M}(x)\ll
x^{1/2+\epsilon}$. Moreover, if the Riemann hypothesis is true, then
\begin{equation}\label{main}
\mathfrak{M}(x)=\sum_{0\leq j<J}\sum_{\rho=1/2+i\gamma\atop \vert
\gamma\vert<T} {x^{\rho/2^j}\over
\rho/2^j}c_j(\rho)+o\big(x^{2^{-J}}\big)
\end{equation}
with some non-zero constants $c_j(\rho)$ if all zeros are simple
(otherwise a modified formula holds with $c_j$ being polynomials in
$\log x$ according to the multiplicities of $\rho$), and
\begin{equation}\label{sounda}
\mathfrak{M}(x)\ll x^{1/2}\exp\big((\log x)^{1/2}(\log\log x)^{14}\big);
\end{equation}
furthermore, the line $\sigma=0$ is a natural boundary for $\mathfrak{m}(s)$.
\end{thm}

It follows from (\ref{main}) that $\mathfrak{M}(x)=\Omega(x^{1/2})$.
In view of Theorem \ref{eins} we deduce that the latter bound holds
also unconditionally.

\noindent {\bf Proof.} If the Riemann hypothesis is true, then
\begin{equation}\label{ui}
\zeta(\sigma+it)^{-1}\ll t^\epsilon\qquad\mbox{for}\quad \sigma>{\textstyle{1\over 2}}
\end{equation}
and all positive $\epsilon$ as $t\to\infty$ (see \cite{titch}, \S
14.2);  moreover, for any real interval of length one, there exists
a real number $t$ from this interval such that the latter estimate
holds also for $s={1\over 2}+it$ (see \cite{titch}, \S 14.16).
Incorporating this bound in place of (\ref{schnee}), we get instead
of estimate (\ref{star})
$$
\mathfrak{M}(x)\ll x^cT^{\epsilon-1}+x^{1-\delta}T^\epsilon+{x\log x\over T}
$$
with any $\delta>{1\over 2}$; now choosing $c$ as in the previous
proof and $T$  such that $T=x^{1/2+\epsilon}$ the desired bound
follows.

If $\mathfrak{M}(x)\ll x^{1/2+\epsilon}$, then we deduce from (\ref{zuz}) that
$$
\mathfrak{m}(s)=\sum_{n=1}^\infty{\mu_{\infty}(n)\over
n^s}=s\int_1^\infty  \mathfrak{M}(u)u^{-s-1}\d u
$$
is convergent and hence analytic for $\sigma>1$, which implies the
non-vanishing of the zeta-function in this half-plane.

For the sake of simplicity, we assume besides the Riemann hypothesis
that all $\zeta$-zeros are simple. Similarly to the proof of Theorem
\ref{eins}, by the calculus of residues, for $\delta:=3\cdot
2^{-J-2}$,
$$
\mathfrak{M}(x)={1\over 2\pi i}\left\{\int_{c-iT}^{\delta-iT}+
\int_{\delta-iT}^{\delta+iT}+\int_{\delta+iT}^{c+iT}\right\}\mathfrak{m}(s){x^s\over
s}\d s+{\mathcal E}+\Sigma,
$$
where ${\mathcal E}$ is the error term bounded in (\ref{e}),
$\Sigma$ is the sum of residues, and the parameter $T$ is chosen
such that $T\neq 2^{-j}\gamma$ for all ordinates $\gamma$ of
$\zeta$-zeros and $0\leq j<J$. All residues arise from zeros
$\rho={1\over 2}+i\gamma$ of $\zeta(s)$; hence
\begin{eqnarray*}
\lefteqn{{\rm Res}_{s=2^{-j}\rho}\mathfrak{m}(s){x^s\over s}}\\
&=&\lim_{s\to
2^{-j}\rho}(s-2^{-j}\rho)\zeta(2^{-j}s)^{-1}N_J(s)\prod_{0\leq\iota<J\atop
\iota\neq j}\zeta(2^{-\iota}s)^{-1}{x^s\over
s}=c_j(\gamma){x^{2^{-j}\rho}\over 2^{-j}\rho}
\end{eqnarray*}
with some non-zero constant $c_j(\gamma)$. Summing up over all zeros
$\rho={1\over 2}+ i\gamma$ with $\vert\gamma\vert$ for all $0\leq
j<J$ yields an expression for $\Sigma$ which constitutes the main
term of the formula (\ref{main}). In order to bound the integrals we
recall that $N_J(s)\ll 1$ for $\sigma>2^{-J}$ (hence on
$\sigma=\delta$) by (\ref{tannenbaum}). Under assumption of the
Riemann hypothesis we have, besides (\ref{ui}),
$$
\zeta(\sigma+it)^{-1}\ll t^{\sigma-1/2+\epsilon} \qquad \mbox{for}\quad \sigma<{\textstyle{1\over 2}};
$$
this follows easily from (\ref{ui}) by use of the functional equation and Stirling's formula; hence
$$
\prod_{0\leq j<J}\zeta(2^js)^{-1}\ll t^{(2^J-1)\sigma-J/2+\epsilon}
$$
as $t\to+\infty$. In view of (\ref{mn}) we find
\begin{eqnarray*}
\lefteqn{\int_{\delta-iT}^{\delta+iT}\mathfrak{m}(s){x^s\over s}\d s}\\
&=&\int_{\delta-iT}^{\delta+iT}N_{J+1}(s)\zeta(2^Js)^{-1}\prod_{0\leq
j<J} \zeta(2^js)^{-1}{x^s\over s}\d s\ll x^\delta
T^{(2^J-1)\delta-J/2+\epsilon}
\end{eqnarray*}
and
$$
\int_{\delta\pm iT}^{c\pm iT}\mathfrak{m}(s){x^s\over s}\d s\ll x T^{\epsilon-1}.
$$

Collecting together and chosing $c$ as in the previous proof, we arrive at
$$
\mathfrak{M}(x)-\Sigma\ll {x\log x\over T}+xT^{\epsilon-1}+x^\delta T^{(2^J-1)\delta-J/2+\epsilon}.
$$
Since $(2^J-1)\delta-J/2<0$, the right-hand side is $=o(x^{2^{-J}})$ which proves (\ref{main}).
\par

Estimate (\ref{sounda}) is a consequence of recent work of Soundararajan \cite{sound}
who obtained for the summatory function of the ordinary M\"obius function $\mu(n)$ via Perron's formula
$$
\sum_{n\leq x}\mu(n)={1\over 2\pi i}\int_{c-i[x]}^{c+i[x]}{x^s\over s\zeta(s)}\d s+O(\log x)
$$
the estimate on the right-hand side of (\ref{sounda}) by contour
integration to the right of the critical line $\sigma={1\over 2}$;
since $\mathfrak{M}(x)$ is also bounded by the above integral for
$\sigma>{1\over 2}$, we may adopt his bound for our case too.

For the assertion that there is no meromorphic continuation beyond the imaginary
axis it suffices to show that in any neighbourhood of any point $it$ with large
imaginary part $t>0$ there exists a pole of $\mathfrak{m}(s)$. Given any
$\epsilon>0$, we have to find a nontrivial zero $\rho=\beta+i\gamma$ of $\zeta(s)$ such that
$$
\vert it-2^{-j}(\beta+i\gamma)\vert<\epsilon
$$
for some positive integer $j$. Since for any sufficiently large $j$ we have
$0<2^{-j}\beta\leq {\epsilon\over 2}$, we have to find a zero $\rho=\beta+i\gamma$ satisfying
\begin{equation}\label{condi}
\vert 2^jt-\gamma\vert <2^{j-1}\epsilon.
\end{equation}

By the Riemann--von Mangoldt formula with an error term under
assumption of the truth of the Riemann hypothesis,
$$
N(T)={T\over 2\pi}\log{T\over 2\pi e}+O\left({\log T\over\log\log T}\right)
$$
(see \cite{titch}, \S 14.13), we find for the number of zeros satisfying condition (\ref{condi})
the estimate
$$
N(2^jt+2^{j-1}\epsilon)-N(2^jt-2^{j-1}\epsilon)\geq {2^j\epsilon
\over 2\pi} \log{2^jt-2^{j-1}\epsilon\over 2\pi e}+O\left({\log
(2^jt)\over \log\log (2^jt)}\right).
$$

Setting $2^{1-j}=\epsilon$, this leads to
$$
N(2t/\epsilon+1)-N(2t/\epsilon-1)\geq {1\over \pi}\log(2t/\epsilon)+o(\log (2t/\epsilon)),
$$
which is positive for sufficiently large $t$. Hence, there is a
singularity in any neighbourhood of almost any arbitrary point $it$,
and thus the imaginary axis is a natural boundary for
$\mathfrak{m}(s)$. The theorem is proved. $\bullet$

Formula (\ref{main}) is similar to the explicit formula for the summatory function
of the classical M\"obius $\mu$-function:
$$
{\sf M}(x)=\sum_{n\leq x}\mu(n)=\sum_{0<\gamma<T} {x^\rho\over \rho\zeta'(\rho)}+{\rm error}(x,T),
$$
which is valid under assumption of the Riemann hypothesis and the so-called
essential simplicity hypothesis that all $\zeta$-zeros are simple, resp.
with obvious modifications if there are multiple zeros (see \cite{titch}, \S 14.27).

\section{Heuristics}

Finally, we discuss some related heuristics based on an old idea due to Denjoy \cite{denj} for the classical $\mu$-function which give support for Riemann's hypothesis. Whereas Denjoy argued for the M\"obius $\mu$-function we consider the modified function $\mu_{\infty}(n)$. Assume that $\{X_n\}$ is a sequence of random variables with distribution
$$
{\bf P}(X_n=+1)={\bf P}(X_n=-1)={\textstyle{1\over 2}}.
$$

Define $S_0=0$ and $S_n=\sum_{j=1}^nX_j$, then $\{S_n\}$ is a
symmetrical random walk in $\Z^2$ with starting point at $0$. A
simple application of Chebyshev's inequality yields, for any
positive $c$,
$$
{\bf P}\{\vert S_n\vert\geq c n^{1/2}\}\leq {1\over 2c^2},
$$
which shows that {\it large} values for $S_n$ are {\it rare} events.
By the theorem of Moivre-Laplace this can be made more precise. It
follows that
$$
\lim_{n\to\infty}{\bf P}\left\{\vert S_n\vert< cn^{1/2}\right\}=
{1\over\sqrt{2\pi}}\int_{-c}^c\exp\left(-{\textstyle{1\over
2}}x^2\right)\d x.
$$

Since the right-hand side above tends to $1$ as $c\to\infty$, we
obtain
$$
\lim_{n\to\infty}{\bf P}\left\{\vert S_n\vert\ll n^{1/2+\epsilon}\right\}=1
$$
for every $\epsilon>0$. We observe that this might be regarded as a
model  for the value-distribution of the modified M\"obius function
$\mu_{\infty}(n)$. (However, for the classical M\"obius function
$\mu(n)$ one has to exclude the squarefull integers $n$ since for
those values $\mu(n)=0$.) The law of the iterated logarithm would
even give the stronger estimate
$$
\lim_{n\to\infty}{\bf P}\left\{\vert S_n\vert\ll (n\log\log n)^{1/2}\right\}=1,
$$
which suggests for $\mathfrak{M}(x)$ the upper bound $(x\log\log
x)^{1/2}$. This estimate is pretty close to the
$\mu_{\infty}$-variant of the so-called weak Mertens hypothesis:
$$
\int_1^X \left({\mathfrak{M}(x)\over x}\right)^2\d x\ll \log X.
$$

The latter bound implies the Riemann hypothesis and the essential simplicity
hypothesis; the proof follows exactly the same argument as in the case of the
classical M\"obius function (see \cite{titch}, \S 14.29) since the generating
Dirichlet series of $\mu$ and $\mu_{\infty}$ differ only by the factor
$\prod_{j\geq 1}\zeta(2^js)^{-1}$ which has no deeper influence.

\end{document}